\documentclass[11pt]{amsart}

\usepackage{graphicx}
\usepackage{epsfig}
\usepackage{amsmath}
\usepackage{amssymb}
\usepackage{bigints}
\usepackage{tikz}
\usepackage{graphicx}
\usepackage{epsfig}
\usepackage[noadjust]{cite}
\usepackage{xcolor}

\usepackage[colorlinks,citecolor=red,urlcolor=blue,bookmarks=false,hypertexnames=true]{hyperref}

\theoremstyle{plain}
\newtheorem{theorem}                {Theorem}      [section]
\newtheorem*{theorem*}                {Theorem \ref{thm:appl}}

\newtheorem{corollary}    [theorem]  {Corollary}
\newtheorem{lemma}        [theorem]  {Lemma}

\theoremstyle{definition}

\newtheorem{remark}       [theorem]  {Remark}

\DeclareMathOperator{\trace}{trace} 

\DeclareMathOperator{\ricci}{Ricci}

\numberwithin{equation}{section}

\begin{document}

\title[A new gap for $CMC$ biharmonic hypersurfaces in Euclidean spheres]
{A new gap for $CMC$ biharmonic hypersurfaces in Euclidean spheres}

\author{Simona~Nistor}

\address{Faculty of Mathematics\\ Al. I. Cuza University of Iasi\\
Blvd. Carol I, 11 \\ 700506 Iasi, Romania} \email{nistor.simona@ymail.com}

\thanks{This work was supported by a grant of the Romanian Ministry of Research and Innovation, CNCS – UEFISCDI, project number PN-III-P1-1.1-PD-2019-0429, within PNCDI III}



\begin{abstract}
In this note we improve a gap result concerning the range of the mean curvature of complete $CMC$ proper-biharmonic hypersurfaces in unit Euclidean spheres.
\end{abstract}

\thanks {I would like to thank professors Cezar Oniciuc and Yu Fu for useful discussions and suggestions.}

\maketitle

\section{Introduction}

The study of biharmonic maps is of great interest for many mathematicians, especially for those working in Differential Geometry. As suggested by J. Eells and J. H. Sampson (see \cite{ES64,ES65}), or J. Eells and L. Lemaire (see \cite{EL83}), a biharmonic map $\varphi:M\to N$ between two Riemannian manifolds is a critical point of the bienergy functional
$$
E_2:C^\infty(M,N)\to \mathbb{R}, \qquad E_2(\varphi)=\frac{1}{2}\int_{M}\left|\tau(\varphi)\right|^2 \ dv,
$$
where $M$ is compact and $\tau(\varphi)=\trace \nabla d\varphi$ is the tension field associated to $\varphi$. In 1986, G.-Y. Jiang (see \cite{J86,J86-2}) proved that the biharmonic maps are characterized by the vanishing of their bitension fields, where the bitension field associated to a map $\varphi$ is given by
$$
\tau_2(\varphi)=-\Delta \tau(\varphi)-\trace R^N\left(d\varphi(\cdot),\tau(\varphi)\right)d\varphi(\cdot).
$$
The nonlinear fourth order elliptic equation $\tau_2(\varphi)=0$ is called the biharmonic equation.

Trivially, any harmonic map is biharmonic, so we will focus on the study of proper-biharmonic maps, i.e., the biharmonic maps which are not harmonic. A biharmonic submanifold of $N^n$ is an isometric immersion $\varphi:M^m\to N^n$ which is also a biharmonic map. Sometimes, throughout this paper, we will indicate a submanifold as $M^m$ rather than mentioning the isometric immersion $\varphi$.

When the curvature of the ambient manifold is non-positive, with just one exception (see \cite{OT12}), we have only non-existence results, i.e., biharmonicity implies harmonicity (minimality). On the other hand, in spaces of positive curvature, especially in Euclidean spheres, we have many examples and classification results for proper-biharmonic submanifolds (see, for example \cite{BMO12, FO20, FYZ22, OHTh12, OC20}).

A particular case is given by the study of biharmonic hypersurfaces in the unit Euclidean sphere $\mathbb{S}^{m+1}$. A very interesting situation is the case when the mean curvature $f$ is a non-zero constant because, with this assumption, the hypersurface $M^m$ is proper-biharmonic if and only if $|A|^2=m$, where $A$ is the shape operator. Since the minimal, i.e., $f=0$, hypersurfaces with $|A|^2=m$ were already classified by S. S. Chern, M. do Carmo and S. Kobayashi in their famous paper \cite{CCK70}, the study of proper-biharmonic hypersurfaces in Euclidean spheres with $f$ constant, i.e., $CMC$, can be seen as a natural generalization of the above mentioned classical problem.

The only known examples of proper-biharmonic hypersurfaces in $\mathbb{S}^{m+1}$ are open parts of the small hypersphere of radius $1/\sqrt{2}$, i.e.,  $\mathbb{S}^m\left(1/\sqrt{2}\right)$, and of the Clifford tori $S^{m_1}\left(1/\sqrt{2}\right)\times S^{m_2}\left(1/\sqrt{2}\right)$, with $m_1\neq m_2$ and $m_1+m_2=m$ (see \cite{CMO01, J86-2}).

Moreover, it was proved that, under various additional geometric assumptions, the proper-biharmonic hypersurfaces have to be the above ones or (at least) they must be $CMC$. Consequently, the following two conjectures have been proposed in $2008$ (see \cite{BMO08}).

\textit{\textbf{Conjecture 1. (C1)} Any proper-biharmonic hypersurface in $\mathbb{S}^{m+1}$ is either an open part of $\mathbb{S}^m\left(1/\sqrt{2}\right)$, or an open part of $S^{m_1}\left(1/\sqrt{2}\right)\times S^{m_2}\left(1/\sqrt{2}\right)$, with $m_1\neq m_2$ and $m_1+m_2=m$.}

\textit{\textbf{Conjecture 2. (C2)} Any proper-biharmonic submanifold in $\mathbb{S}^{m+1}$ is $CMC$.}

We recall that, when the ambient space is the Euclidean space $\mathbb{R}^n$, the famous \textit{Chen Conjecture} remains unsolved. The conjecture says: \textit{any biharmonic submanifold in $\mathbb{R}^n$ is minimal} (see \cite{C91}). Since any $CMC$ biharmonic submanifold in $\mathbb{R}^n$ is minimal (see \cite{D92}), the \textit{Chen Conjecture} can be reformulated in a weaker form:

\textit{\textbf{Weak version of Chen Conjecture.} Any biharmonic submanifold in $\mathbb{R}^n$ is $CMC$.}

Thus, conjecture \textit{C2} can be seen as an extension of the \textit{Weak version of Chen Conjecture} to Euclidean spheres.

Obviously, \textit{C2} (for hypersurfaces) is weaker than \textit{C1}, but to directly prove the first conjecture seems to be quite a complicated task. Conjecture \textit{C2} can be seen as an intermediate step for proving \textit{C1}. However, even if \textit{C2} will be proved, the proof of \textit{C1} will still be a real challenge, and the additional hypothesis of compactness of the hypersurface does not simplify it.

Until now, \textit{C1} was proved only when $m=2$ (see \cite{CMO01}). Also, it was proved in several particular cases, imposing additional geometric hypotheses on the biharmonic hypersurface: $m=3$ and the hypersurface $M$ is complete (see \cite{BMO10, OHTh12}), or $M^m$ has at most two distinct principal curvatures at any point (see \cite{BMO08}), or $M^m$ is isoparametric (see \cite{IIU09,IIU10}), or $M^m$ is $CMC$ and has non-positive sectional curvature (see \cite{OHTh12}), or $M^m$ is compact and belongs to a hemisphere (see \cite{V20}), etc.

As it is not clear how \textit{C2} could imply \textit{C1}, we propose an intermediary objective.

\textit{\textbf{Open Problem.} Let $M^m$ be a $CMC$ proper-biharmonic hypersurface in $\mathbb{S}^{m+1}$. Then, the set of all possible values of the mean curvature is discreet and, more precisely,
\begin{equation*}
f\in\left\{\frac{m-2r}{m} \quad | \quad r\in\mathbb{N},\ 0\leq r\leq s^\ast\right\},
\end{equation*}
where $s^\ast=s-1$, if $m=2s$, and $s^\ast=s$, if $m=2s+1$.}

Finally, we should prove that, if $f=(m-2r)/m$, then $M^m$ must be an open part of $S^{r}\left(1/\sqrt{2}\right)\times S^{m-r}\left(1/\sqrt{2}\right)$.

We mention here that there is a deep link between the proof of \textit{C1} knowing that \textit{C2} is true and the \textit{Generalized Chern Conjecture}, as we will explain below.

First, inspired by the well-known \textit{Chern Conjecture} concerning compact minimal hypersurfaces and its generalization to $CMC$ hypersurfaces (see, for example \cite{dABSW20, LXX21, SWY12}), we can state the following

\textit{\textbf{Generalized Chern Conjecture.} Let $M^m$ be a $CMC$ hypersurface in $\mathbb{S}^{m+1}$ with constant squared norm of the shape operator. Then, $M$ is isoparametric.}

Further, as a non-minimal $CMC$ hypersurfaces in $\mathbb{S}^{m+1}$ is proper-biharmonic if and only if $|A|^2=m$, if the \textit{Generalized Chern Conjecture} and \textit{C2} will be proved to be true, then our \textit{C1} will follow immediately using the results in \cite{IIU09, IIU10}. Therefore, the proof of \textit{C1} under the $CMC$ assumption can be seen as a special case of \textit{Generalized Chern Conjecture}.

However, the \textit{Generalized Chern Conjecture} seems very difficult to be proved in its full generality, so we think that there are more chances to prove it under the additional hypothesis $|A|^2=m$, i.e. to prove \textit{C1}, assuming that \textit{C2} is true.

The main result of this paper gives a partial answer to the above \textit{Open Problem} showing that
$$
f\in \left(0,\gamma\right)\cup \left\{\frac{m-2}{m}\right\}\cup \{1\},
$$
where $\gamma$ is a real constant depending only on $m$ and
$$
\frac{m-3}{m}<\gamma<\frac{m-2}{m}.
$$

\textbf{Conventions and notations.} In this paper, the Riemannian metrics are indicated by $\langle\cdot,\cdot\rangle$. We assume that all manifolds are connected and oriented, and we use the following sign conventions for the rough Laplacian acting on sections of $\varphi^{-1}(TN)$ and for the curvature tensor field, respectively:
$$
\Delta^{\varphi}=-\trace\left(\nabla^{\varphi}\nabla^{\varphi}-\nabla^{\varphi}_{\nabla}\right)
$$
and
$$
R(X,Y)Z=[\nabla_X,\nabla_Y]Z-\nabla_{[X,Y]}Z.
$$
For an (oriented) hypersurface $\varphi:M^m\to N^{m+1}$ we label the principal curvatures of $M$ such that
$$
\lambda_1\leq\lambda_2\leq\cdots\leq\lambda_m.
$$

\section{Proper-biharmonic hypersurfaces with constant mean curvature}

Let $M^m$ be a hypersurface in the unit Euclidean sphere $\mathbb{S}^{m+1}$. For simplicity, we assume that $M$ is oriented. We recall that $A=A_\eta$ denotes the shape operator of $M$ and $B$ is the associated second fundamental form, $B(X,Y)=\langle A(X),Y\rangle\eta$, where $\eta$ is the unit normal vector field, globally defined on $M$. In this paper, we consider the normalized version for the mean curvature function, i.e., $f=\left(\trace A\right)/m$. When $M$ is a non-minimal $CMC$ hypersurface, i.e., $f$ is a non-zero constant, we can consider $\eta=H/|H|$ and so $f=|H|$, where $H$ is the mean curvature vector field.

Now, we recall a first result that supports the above \textit{Open Problem}.

\begin{theorem}[\cite{ODocTh03, OHTh12}]\label{th:0-1}
Let $\varphi:M^m\to\mathbb{S}^n$ be a $CMC$ proper-biharmonic submanifold. Then $f\in (0,1]$, and $f=1$ if and only if $\varphi$ induces a minimal immersion of $M$ into $\mathbb{S}^{n-1}\left(1/\sqrt{2}\right)\subset\mathbb{S}^{n}$.
\end{theorem}

We also recall

\begin{theorem}[\cite{BO12}]\label{th:gap1}
Let $\varphi:M^m\to\mathbb{S}^n$ be a proper-biharmonic immersion with parallel mean curvature vector field. Assume that $m>2$ and $f\in(0,1)$. Then $f\in\left(0,(m-2)/m\right]$, and $f=(m-2)/m$ if and only if locally $\varphi(M)$ is an open part of a standard product
$$
\mathbb{S}^1\left(\frac{1}{\sqrt{2}}\right)\times M_1\subset\mathbb{S}^n,
$$
where $M_1$ is a minimal embedded submanifold of $\mathbb{S}^{n-2}\left(1/\sqrt{2}\right)$. Moreover, if $M$ is complete, then the above decomposition of $\varphi(M)$ holds globally, where $M_1$ is a complete minimal submanifold of $\mathbb{S}^{n-2}\left(1/\sqrt{2}\right)$.
\end{theorem}

As a corollary of Theorem \ref{th:gap1}, we get a second result confirming the \textit{Open Problem}.

\begin{corollary}[\cite{BO12}]\label{cor:gapCMC}
Let $\varphi:M^m\to\mathbb{S}^{m+1}$ be a $CMC$ proper-biharmonic hypersurface with $m>2$. Then $f\in\left(0,(m-2)/m\right]\cup\left\{1\right\}$. Moreover, $f=1$ if and only if $\varphi(M)$ is an open subset of the small hypersphere $\mathbb{S}^m\left(1/\sqrt{2}\right)$, and $f=(m-2)/m$ if and only if $\varphi(M)$ is an open subset of the standard product $\mathbb{S}^1\left(1/\sqrt{2}\right)\times\mathbb{S}^{m-1}\left(1/\sqrt{2}\right)$.
\end{corollary}

The same result of Theorem \ref{th:gap1} was proved, independently and in the same time, in \cite{WW12}.

We mention that Corollary \ref{cor:gapCMC} can be reobtained, under the additional assumption $M^m$ compact, from the following known result.

\begin{theorem}[\cite{CN90,X90}]\label{th:A}
Let $M^m$ be a compact non-minimal $CMC$ hypersurface in $\mathbb{S}^{m+1}$. If $|A|^2\leq\alpha$, then $M$ is either an umbilical sphere, or a Clifford torus.
\end{theorem}

Here,
\begin{equation}\label{eq:alpha}
\alpha=\alpha (m,f)=m+\frac{m^3}{2(m-1)}f^2-\frac{m(m-2)}{2(m-1)}\sqrt{m^2f^4+4(m-1)f^2}.
\end{equation}

Indeed, if $M^m$ is a compact $CMC$ proper-biharmonic hypersurface in $\mathbb{S}^{m+1}$, i.e., $|A|^2=m$, the hypothesis $|A|^2\leq\alpha$ in the above theorem is equivalent to $f\geq(m-2)/m$. Thus, taking into account Theorem \ref{th:0-1}, we get that a compact $CMC$ proper-biharmonic hypersurface satisfying $|A|^2\leq\alpha$ must obey $f\in\left[(m-2)/m,1\right]$. Therefore, when $f$ belongs to the above interval, from Theorem \ref{th:A} we get that $M$ is either an umbilical sphere, or a Clifford torus. It is known that the only biharmonic umbilical hypersurface is $\mathbb{S}^m\left(1/\sqrt{2}\right)$, and so $f=1$, and the only biharmonic Clifford torus is $\mathbb{S}^{m_1}\left(1/\sqrt{2}\right)\times \mathbb{S}^{m_2}\left(1/\sqrt{2}\right)$, where $1\leq m_1\leq m-1$, $m_1\neq m_2$ and $m_1+m_2=m$ (see \cite{BMO10,J86-2}). In the later case, we can assume $m_1<m_2$ and get $f=\left(m_2-m_1\right)/m$. Further, as by hypothesis $f$ belongs to the interval $\left[(m-2)/m,1\right]$, the only possibility is that $m_1=1$, and therefore $f=(m-2)/m$.

Our main result is an improvement of Corollary \ref{cor:gapCMC}. We show that there is a larger gap for $f$ than $\left((m-2)/m,1\right)$. More precisely, considering $m\geq 4$ and denoting
\begin{equation}\label{gamma}
\gamma=(m-2)\sqrt{\frac{m-1}{m^2(m-1)+B_m\left(B_m+m^2\right)}},
\end{equation}
where
\begin{equation}\label{Bm}
B_m=\left\{
\begin{array}{ll}
  0.2, & 4 \leq  m \leq 42 \\\\
  0.199, &  43 \leq m \leq 65\\\\
  0.198, &  66 \leq m \leq 149\\\\
  0.197, &   m \geq 150
\end{array}
\right.,
\end{equation}
we have

\begin{theorem}\label{th:mainTh}
Let $M^m$ be a complete $CMC$ proper-biharmonic hypersurface in $\mathbb{S}^{m+1}$. If $m\geq 4$ and the mean curvature $f\in \left[\gamma,(m-2)/m\right]$, then $f=(m-2)/m$ and $M=\mathbb{S}^{1}\left(1/\sqrt{2}\right)\times \mathbb{S}^{m-1}\left(1/\sqrt{2}\right)$.
\end{theorem}

\begin{remark}
The real constant $\gamma$ satisfies
$$
\frac{m-3}{m}<\gamma<\frac{m-2}{m},\qquad m \geq 4.
$$
\end{remark}

A direct consequence of Theorem \ref{th:mainTh} is the next result, which gives the new gap for $f$.

\begin{corollary}\label{cor:mainCor}
Let $M^m$ be a complete $CMC$ proper-biharmonic hypersurface in $\mathbb{S}^{m+1}$, with $m\geq 4$. Then
$$
f\in(0,\gamma)\cup\left\{\frac{m-2}{m}\right\}\cup\{1\}.
$$
Moreover, $f=(m-2)/m$ if and only if $M=\mathbb{S}^{1}\left(1/\sqrt{2}\right)\times \mathbb{S}^{m-1}\left(1/\sqrt{2}\right)$, and $f=1$ if and only if $M=\mathbb{S}^{m}\left(1/\sqrt{2}\right)$.
\end{corollary}

We note that Theorem \ref{th:mainTh} can be deduced from Theorem 1 in \cite{GLX18} taking $|A|^2=m$, $f\in\left[\gamma,(m-2)/m\right]$, but more precise values of $B_m$. We will provide a slightly simpler proof than the original one, using the properties of biharmonic hypersurfaces.

\section{The proof of Theorem \ref{th:mainTh}}

In the first part of the proof we will give some algebraic results that will be very helpful. They will mainly involve certain quantities denoted by $\mu_i$, with $ i\in\{1,2,\ldots, m\}$, $\phi$, $\eta$ and $\sigma$. At the beginning, we will justify why we introduce them.

We recall that when $M=\mathbb{S}^{1}\left(1/\sqrt{2}\right)\times \mathbb{S}^{m-1}\left(1/\sqrt{2}\right)$, the principal curvatures are constant and given by
$$
\lambda_1=-1 \quad \text{ and } \quad  \lambda_2=\lambda_3=\cdots=\lambda_m=1=\frac{mf-\lambda_1}{m-1}.
$$

For a hypersurface $M^m$ in $\mathbb{S}^{m+1}$ satisfying the hypotheses of Theorem \ref{th:mainTh}, at a certain step, the term
$$
\trace A^3=\sum_{i=1}^{m}\lambda_i^3
$$
will appear. For our objective, it will be more convenient to replace $\lambda_i$ by
$$
\lambda_i= \sqrt{m(1-f^2)}\mu_i+f.
$$
Now, the advantage of using $\mu_i$'s is that
\begin{equation}\label{eqmu}
\sum_{i=1}^{m}\mu_i=0 \quad \text{ and } \quad \sum_{i=1}^{m}\mu_i^2=1.
\end{equation}

We note that, in the particular case when $M=\mathbb{S}^{1}\left(1/\sqrt{2}\right)\times \mathbb{S}^{m-1}\left(1/\sqrt{2}\right)$, we have
$$
\mu_1=-\sqrt{\frac{m-1}{m}}\quad \text{ and } \quad \mu_2=\mu_3=\cdots=\mu_m=\frac{1}{\sqrt{m(m-1)}}=-\frac{\mu_1}{m-1}.
$$

We also need to consider on the hypersurface $M^m$ three functions $\phi$, $\eta$ and $\sigma$ given by
\begin{equation}\label{eq:PhiEtaSigma}
\phi=\sum_{i=1}^{m}\mu_i^3+\frac{m-2}{\sqrt{m(m-1)}},\quad \eta=\sqrt{\frac{m}{m-1}}\mu_1+1, \quad \sigma=\sqrt{\sum_{i=2}^{m}\left(\mu_i+\frac{\mu_1}{m-1}\right)^2}.
\end{equation}

\begin{remark}
According to Okumura Lemma (see \cite{O74}) we have
$$
0 \leq \phi \leq \frac{2(m-2)}{\sqrt{m(m-1)}}.
$$
As we will see in the following, $\phi$ is indeed nonnegative (Lemma \ref{lemma1}) and in Lemma \ref{lemma2} we will impose an upper bound of $\phi$ which is less than $2(m-2)/\sqrt{m(m-1)}$.
\end{remark}

We mention that, in the particular case, when $M=\mathbb{S}^{1}\left(1/\sqrt{2}\right)\times \mathbb{S}^{m-1}\left(1/\sqrt{2}\right)$, we get
$$
\phi=\eta=\sigma=0.
$$

Fixing arbitrarily a point of a hypersurface $M^m$ in $\mathbb{S}^{m+1}$ that satisfies the hypotheses of Theorem \ref{th:mainTh}, the above functions become, obviously, real numbers, and in the following we will give three algebraic lemmas. The first two lemmas  are obtained in \cite{GLX18}. Their proofs are elementary but skilful and we do not present them here. Concerning the third lemma, we mention that it is originated in a result in \cite{GLX18}, but our statement is  more accurate and the proof is different.

\begin{lemma}[\cite{GLX18}]\label{lemma1}
If $m\geq3$, the real numbers $\phi$, $\eta$ and $\sigma$ satisfy
\begin{equation}\label{eq1Lemma1}
\frac{\sqrt{m(m-1)}}{m-2}\phi\geq\eta\geq\frac{\sigma^2}{2}
\end{equation}
and
\begin{equation}\label{eq2Lemma1}
\phi\sqrt{m(m-1)}\geq\eta\left[3m-3(m+1)\eta-2\sigma\sqrt{m(m-1)}\right].
\end{equation}
\end{lemma}

In order to state the next lemma, we consider a positive number defined by
\begin{eqnarray}\label{eq:alpha0}
\nonumber % Remove numbering (before each equation)
\alpha_0 &=& \alpha-mf^2= \\
&=& m+\frac{m^3}{2(m-1)}f^2-\frac{m(m-2)}{2(m-1)}\sqrt{m^2f^4+4(m-1)f^2}-mf^2.
\end{eqnarray}

\begin{lemma}[\cite{GLX18}]\label{lemma3}
If $m\geq 3$, then the following equality holds
\begin{equation}\label{eqLemma3}
(m-2)f\sqrt{\frac{m}{m-1}\alpha_0}=m\left(f^2+1\right)-\alpha_0.
\end{equation}
\end{lemma}

From the definitions of $\sigma$ and $\eta$ we get a link between the difference $\mu_2-\mu_1$ and a quantity which contains $\sigma$, $\eta$ and $m$. This we will be useful to prove the third lemma.

\begin{equation}\label{ineq:mu2mu1}
\mu_2-\mu_1\geq\left(1-\eta-\sigma\sqrt{\frac{m-1}{m}}\right)\sqrt{\frac{m}{m-1}}.
\end{equation}

Now, using the real constant $B_m$ given in \eqref{Bm}, which, clearly, depends on $m$, we can state the following lemma.

\begin{lemma}\label{lemma2}
Let $m\geq 4$. If
\begin{equation}\label{eq1Lemma2}
\phi\leq\frac{B_m}{2}\sqrt{\frac{m}{m-1}},
\end{equation}
then $2\sigma+3\eta<3/4$ and
\begin{equation}\label{eq2Lemma2}
\mu_2-\mu_1>\frac{2}{3-3^{-10}}\sqrt{\frac{m}{m-1}}.
\end{equation}
\end{lemma}

\begin{proof}
First, using \eqref{eq1Lemma1}, \eqref{eq1Lemma2} and the hypothesis $m\geq 4$, we get that $\eta\leq B_m$.

Then, using \eqref{eq2Lemma1} and \eqref{eq1Lemma2}, we obtain
$$
3\eta-3\left(1+\frac{1}{m}\right)\eta^2-2\eta\sigma\sqrt{1-\frac{1}{m}} -\frac{B_m}{2}\leq 0.
$$
Further, from \eqref{eqmu} and \eqref{eq:PhiEtaSigma}, it is easy to see that $\sigma=\sqrt{2\eta-\eta^2}$, and $\sigma$, $\eta\in[0,1)$. Previously, we have seen that $\eta\leq B_m$, so, in fact $\eta\in \left[0,B_m\right]$ and $\sigma\in \left[0,\sqrt{2B_m-B_m^2}\right]$.

Therefore, the above inequality can be rewritten as
\begin{equation}\label{ineq:fundam1}
3\eta-3\left(1+\frac{1}{m}\right)\eta^2-2\eta\sqrt{1-\frac{1}{m}}\sqrt{2\eta-\eta^2}-\frac{B_m}{2}\leq 0.
\end{equation}
As $m>0$, it follows that $-2\sqrt{1-1/m}>-2$, and then, from \eqref{ineq:fundam1}, we get
\begin{equation}\label{ineq:fundam2}
3\eta-3\left(1+\frac{1}{m}\right)\eta^2-2\eta\sqrt{2\eta-\eta^2}-\frac{B_m}{2}< 0.
\end{equation}
Further, we continue with an argument which allows us to motivate the values of $B_m$ given in \eqref{Bm}.

Since $m\geq 4$, it follows that $-3\left(1+1/m\right)\geq -15/4$, and from \eqref{ineq:fundam2} we obtain
\begin{equation}\label{ineq:fundam3}
3\eta-\frac{15}{4}\eta^2-2\eta\sqrt{2\eta-\eta^2}-\frac{B_m}{2}< 0.
\end{equation}

Let $B_m=0.2$. Now, as $\eta\leq B_m$, $\eta\in \left[0,0.2\right]$ and from \eqref{ineq:fundam3}, using the program Mathematica$^{\footnotesize \text{\textregistered}}$, we get the maximum range of $\eta$, i.e., $\eta\in (0,0.0446008)$. Next, from $\sigma=\sqrt{2\eta-\eta^2}$, we also get the maximum range of $\sigma$, i.e. $\sigma\in (0,0.295317)$.

Clearly, $2\sigma+3\eta<3/4$.

Further, in order to prove \eqref{eq2Lemma2}, taking into account \eqref{ineq:mu2mu1}, it is enough to have
\begin{equation}\label{ineqsuff}
1-\eta-\sigma\sqrt{\frac{m-1}{m}}>\frac{2}{3-3^{-10}}.
\end{equation}

From the maximum ranges of $\eta$ and $\sigma$ (which depend on the assumption $m\geq 4$ and the chosen value of $B_m$), using again Mathematica, we can see that \eqref{ineqsuff} holds only when $m\leq 22$.

So, if $m\in[4,22]$ and $B_m=0.2$, inequality \eqref{eq2Lemma2} holds.

We continue the proof, assuming now that $m\geq 23$ and $B_m$ has the same value, i.e. $B_m=0.2$. The argument will be similar, but since the minimum value of $m$ changed, also the maximum range of $\eta$ and $\sigma$ would change.

More precisely, as $m\geq 23$, it follows that $-3\left(1+1/m\right)\geq -72/23$, and from \eqref{ineq:fundam2} we get
\begin{equation}\label{ineq:fundam4}
3\eta-\frac{72}{23}\eta^2-2\eta\sqrt{2\eta-\eta^2}-\frac{B_m}{2}< 0.
\end{equation}
Since $B_m=0.2$, using Mathematica, from \eqref{ineq:fundam4}, we obtain $\eta\in (0,0.043933)$ and $\sigma\in (0,0.29315)$.

Clearly, $\eta$ and $\sigma$ belonging to these new intervals also satisfy $2\sigma+3\eta<3/4$ and imposing again \eqref{ineqsuff}, it follows, this time, that $m\leq 39$.

So, if $m\in[23,39]$ and $B_m=0.2$, we also proved that inequality \eqref{eq2Lemma2} holds.

The process will continue in the same way. Considering $m\geq 40$ and $B_m=0.2$, from \eqref{ineq:fundam2}, we obtain
\begin{equation}\label{ineq:fundam5}
3\eta-\frac{123}{40}\eta^2-2\eta\sqrt{2\eta-\eta^2}-\frac{B_m}{2}< 0,
\end{equation}
so $\eta\in (0,0.043875)$ and $\sigma\in(0,0.292962)$. Thus $2\sigma+3\eta<3/4$ and imposing \eqref{ineqsuff}, it follows that $m\leq 42$.

So, if $m\in[40,42]$ and $B_m=0.2$, the inequality \eqref{eq2Lemma2} holds.

Trying to continue the proof with the same argument, we assume now $m\geq 43$ and $B_m=0.2$. In this case, we can see that \eqref{ineqsuff} does not hold, so we cannot motivate that \eqref{eq2Lemma2} holds. For this reason, we need to decrease the value of $B_m$.

Let $m\geq 43$ and $B_m=0.199$. In the same way, first we get that $m\leq 62$ and then, supposing $m\geq 63$ and $B_m=0.199$, we obtain $m\leq 65$. If we continue with $m\geq 66$ and $B_m=0.199$, \eqref{ineqsuff} is not valid and we cannot conclude. Thus, we change again the value of $B_m$.

If we considering $m\geq 66$ and $B_m=0.198$, following the same steps as the above, we can conclude that for $m\leq 149$, $2\sigma+3\eta<3/4$ and the inequality \eqref{eq2Lemma2} holds. Instead, we cannot conclude that \eqref{eq2Lemma2} holds with this value of $B_m$ and $m\geq 150$.

Finally, if we assume $m\geq 150$ and $B_m=0.197$ we obtain that $2\sigma+3\eta<3/4$ and \eqref{ineqsuff} is valid for any $m\geq 150$. Therefore, the conclusion of the lemma is true.

\end{proof}

Our purpose is to show that the functions $\phi$, $\eta$ and $\sigma$ vanish on the hypersurface $\varphi:M^m\to\mathbb{S}^{m+1}$ satisfying the hypotheses of Theorem \ref{th:mainTh}, as they do when $M=\mathbb{S}^{1}\left(1/\sqrt{2}\right)\times \mathbb{S}^{m-1}\left(1/\sqrt{2}\right)$. In order to achieve this, we will prove that the smallest distinct principal curvature has constant multiplicity $1$, so $\lambda_1$ is smooth on $M$. Then, we compute and estimate $\Delta\lambda_1$ and further, using Omori-Yau maximum principle, we conclude that the three functions vanish on $M$. In fact, we will show that $A$ is parallel and from here, using the properties of biharmonic hypersurfaces we obtain that $M=\mathbb{S}^{1}\left(1/\sqrt{2}\right)\times \mathbb{S}^{m-1}\left(1/\sqrt{2}\right)$.

We continue our proof recalling that the principal curvature functions $\lambda_i$'s are continuous on $M$ and the set of the points where the number of distinct principal curvatures is locally constant, denoted by $M_A$, is an open and dense subset in $M$. On a non-empty connected component of $M_A$, which is open in $M_A$ and also in $M$, the number of distinct principal curvatures is constant. Therefore, the multiplicities of distinct principal curvatures are constant and thus, on a connected component of $M_A$, the functions $\lambda_i$'s are smooth and the shape operator $A$ is locally smoothly diagonalizable.

In order to prove that $\lambda_1$ has constant multiplicity $1$, we will employ Lemma \ref{lemma2}, showing that $\mu_1<\mu_2$. For this purpose, we recall the following Simons' type formula, valid on $M^m$, that holds for any $CMC$ hypersurface in $\mathbb{S}^{m+1}$,
\begin{equation}\label{eq:nomizu}
\frac{1}{2}\Delta\left|A\right|^2=-\left|\nabla A\right|^2+m^2f^2+\left|A\right|^2\left(\left|A\right|^2-m\right)-mf\trace A^3.
\end{equation}
Formula \eqref{eq:nomizu} was obtained by K. Nomizu in \cite{NS69}.

Next, we define the positive numbers
$$
m_0=m-mf^2 \quad \text{ and } \quad \delta=\frac{m}{m-1}f^2B_m,
$$
that will be useful for our formulas. We note that the hypotheses of Theorem \ref{th:mainTh} are equivalent to
$$
f\leq \frac{m-2}{m} \Leftrightarrow m\geq\alpha \Leftrightarrow \alpha_0\leq m_0
$$
and
$$
f\geq \gamma \Leftrightarrow m\leq\alpha+\delta \Leftrightarrow m_0\leq\alpha_0+\delta,
$$
where $\alpha$ is defined in \eqref{eq:alpha} and $\alpha_0$ in \eqref{eq:alpha0}.

Now, using the biharmonicity hypothesis, i.e., $|A|^2=m$, from \eqref{eq:nomizu} we obtain, on $M$,
\begin{eqnarray*}
% \nonumber % Remove numbering (before each equation)
\left|\nabla A\right|^2 &=& m^2f^2-mf\trace A^3= \\
&=& m^2f^2-mf\sum_{i=1}^{m}\lambda_i^3,
\end{eqnarray*}
which can be easily rewritten as
$$
\left|\nabla A\right|^2=m_0^2-m\left(f^2+1\right)m_0-mf\sqrt{m_0^3}\sum_{i=1}^{m}\mu_i^3.
$$
Using the definition of $\phi$ and the above expression of $\left|\nabla A\right|^2$, we get, on $M$,
\begin{equation}\label{eq1}
\left|\nabla A\right|^2+mf\phi\sqrt{m_0^3}=m_0\left[m_0-m\left(f^2+1\right)+(m-2)f\sqrt{\frac{m}{m-1}m_0}\ \right].
\end{equation}

In order to fulfill the hypothesis of Lemma \ref{lemma2} and since
\begin{equation}\label{ineq:mfPhiSqrt}
mf\phi\sqrt{m_0^3}\leq \left|\nabla A\right|^2+mf\phi\sqrt{m_0^3},
\end{equation}
the next step consists in finding a convenient upper bound for the term in the right hand side of \eqref{eq1}.

First, as $m_0\leq\alpha_0+\delta$, it is easy to see that $\sqrt{m_0}<\sqrt{\alpha_0}+\delta/\left(2\sqrt{\alpha_0}\right)$, and then, using these inequalities and \eqref{eqLemma3}, we get
\begin{equation}\label{eq2}
m_0-m\left(f^2+1\right)+(m-2)f\sqrt{\frac{m}{m-1}m_0} < \delta+\frac{m-2}{2\sqrt{\alpha_0}}\delta f\sqrt{\frac{m}{m-1}}.
\end{equation}
Clearly, from the definition of $m_0$, we have $m_0<m$ and therefore $\alpha_0<m$. Now, using the definition of $\delta$, the relation \eqref{eqLemma3} and the inequality $\alpha_0-m<0$, we obtain
\begin{eqnarray}\label{ineqTerm1}
\nonumber % Remove numbering (before each equation)
\delta &<& \frac{(m-2)B_m}{m-1}f\sqrt{\frac{m}{m-1}\alpha_0}< \\
&<& B_m f\sqrt{\frac{m}{m-1}\alpha_0},
\end{eqnarray}
and, then
\begin{eqnarray*}
% \nonumber % Remove numbering (before each equation)
\frac{m-2}{2\sqrt{\alpha_0}}\delta f\sqrt{\frac{m}{m-1}} &<& \frac{m-2}{m-1}\frac{B_m}{2}mf^2 < \\
&<& \frac{B_m}{2}mf^2.
\end{eqnarray*}
If we consider again \eqref{eqLemma3} and the inequality $\alpha_0-m<0$, we get
\begin{equation}\label{ineqTerm2}
\frac{m-2}{2\sqrt{\alpha_0}}\delta f\sqrt{\frac{m}{m-1}} < \frac{(m-2)B_m}{2}f\sqrt{\frac{m}{m-1}\alpha_0}.
\end{equation}
Using \eqref{eq2}, \eqref{ineqTerm1}, \eqref{ineqTerm2} and then, $\alpha_0\leq m_0$, one obtains
\begin{eqnarray}\label{ineqF}
\nonumber
m_0-m\left(f^2+1\right)+(m-2)f\sqrt{\frac{m}{m-1}m_0} &<& \frac{m B_m}{2}f\sqrt{\frac{m}{m-1}\alpha_0}\leq \\
&\leq & \frac{m B_m}{2}f\sqrt{\frac{m}{m-1}m_0}.
\end{eqnarray}
We note that, from \eqref{eq1}, \eqref{ineq:mfPhiSqrt} and \eqref{ineqF}, we have on $M$
\begin{eqnarray*}
% \nonumber % Remove numbering (before each equation)
mf\phi\sqrt{m_0^3} &\leq &\left|\nabla A\right|^2+mf\phi\sqrt{m_0^3}< \\
&<&\frac{m B_m}{2}f\sqrt{\frac{m}{m-1}m_0^3},
\end{eqnarray*}
and therefore, on $M$,
$$
\phi<\frac{B_m}{2}\sqrt{\frac{m}{m-1}}.
$$

Now, we can apply Lemma \ref{lemma2} and achieve $\mu_2>\mu_1$, which is equivalent to $\lambda_2>\lambda_1$ on $M^m$. Therefore, since the smallest principal curvature $\lambda_1$ of $M$  has (constant) multiplicity $1$ on $M$, it follows that it is smooth on $M$ and there exists a local smooth unit vector field $E_1$ such that $A\left(E_1\right)=\lambda_1E_1$ (see \cite{N73}). Then, we can find a local expression of $\Delta \lambda_1$ but, in order to work with, it is more convenient to fix arbitrarily a point $p$ and consider $\left\{e_1,\ldots,e_m\right\}$ an orthonormal basis which diagonalize the shape operator $A$ such that $e_1=E_1(p)$. Let $b_{ijk}=\langle \left(\nabla A\right)\left(e_i,e_j\right),e_k\rangle$ be the components of the totally symmetric tensor $\langle \left(\nabla A\right)\left(\cdot,\cdot\right),\cdot\rangle$. It was shown in \cite{GLX18} that, at $p$,
\begin{equation*}
\Delta\lambda_{1}=mf+\left(|A|^2-m\right)\lambda_{1}-mf\lambda_{1}^2-2\sum\limits_{\substack{i=1 \\ k\geq 2}}^m \frac{b_{i1k}^2}{\lambda_{1}-\lambda_{k}}.
\end{equation*}
Now, since $M^m$ is proper-biharmonic and therefore $|A|^2=m$, the above equation (that holds for any $CMC$ hypersurface in $\mathbb{S}^{m+1}$) becomes
\begin{equation}\label{eq:laplaceanLp}
\Delta\lambda_{1}= mf-mf\lambda_{1}^2-2\sum\limits_{\substack{i=1 \\ k\geq 2}}^m \frac{b_{i1k}^2}{\lambda_{1}-\lambda_{k}}.
\end{equation}

Using the link between $\lambda_i$'s and $\mu_i$'s and between $\mu_1$ and $\eta$, we can rewrite the above expression as
$$
\Delta\lambda_1= mf-mf\left[(\eta-1)\sqrt{\frac{m-1}{m}m_0}+f\right]^2-\frac{2}{\sqrt{m_0}}\sum\limits_{\substack{i=1 \\ k\geq 2}}^m \frac{b_{i1k}^2}{\mu_1-\mu_k}.
$$
By a straightforward computation, combining in a suitable way the terms from the right hand side of the above equation, and using the definition of $m_0$, one gets, at $p$
\begin{eqnarray*}
% \nonumber % Remove numbering (before each equation)
\Delta\lambda_1 \quad = & - & \eta\sqrt{m_0}\left[(\eta-2)(m-1)f\sqrt{m_0}+2f^2\sqrt{m(m-1)}\right]- \\
  & - & \left[m_0-m\left(f^2+1\right)+(m-2)f\sqrt{\frac{m}{m-1}m_0}\right]\sqrt{\frac{m-1}{m}m_0}- \\
  & - & \frac{2}{\sqrt{m_0}}\sum\limits_{\substack{i=1 \\ k\geq 2}}^m \frac{b_{i1k}^2}{\mu_1-\mu_k}.
\end{eqnarray*}
We notice that the second squared parenthesis can be replaced by $\left(\left|\nabla A\right|^2+mf\phi\sqrt{m_0^3}\right)/m_0$ from \eqref{eq1}, and we can rewrite the above equation as
\begin{eqnarray}\label{eq:deltaLambda1}
 \nonumber % Remove numbering (before each equation)
\Delta\lambda_1 \quad = && \eta \sqrt{m_0}\left\{-(\eta-2)(m-1)f\sqrt{m_0}+\left[m_0-m\left(f^2+1\right)\right]\sqrt{\frac{m-1}{m}}\right\}- \\
& - & \sqrt{\frac{m-1}{m m_0}}\left|\nabla A\right|^2-m_0\phi f\sqrt{m(m-1)}-\frac{2}{\sqrt{m_0}}\sum\limits_{\substack{i=1 \\ k\geq 2}}^m \frac{b_{i1k}^2}{\mu_1-\mu_k}.
\end{eqnarray}
Further, in order to obtain $\Delta\lambda_1<0$, we will find certain convenient upper bounds for some terms in the right hand side of \eqref{eq:deltaLambda1}.

First, it is easy to see that
\begin{eqnarray*}
% \nonumber % Remove numbering (before each equation)
\Theta &:=& -(\eta-2)(m-1)f\sqrt{m_0}+\left[m_0-m\left(f^2+1\right)\right]\sqrt{\frac{m-1}{m}} <\\
   &<& 2(m-1)f\sqrt{m_0}+\left[m_0-m\left(f^2+1\right)\right]\sqrt{\frac{m-1}{m}}.
\end{eqnarray*}
Second, using $m_0\leq \alpha_0+\delta$, $\sqrt{m_0}<\sqrt{\alpha_0}+\delta/ \left(2\sqrt{\alpha_0}\right)$ and
$$
f^2+1=\frac{1}{m}\left((m-2)f\sqrt{\frac{m}{m-1}\alpha_0}+\alpha_0\right),
$$
obtained from Lemma \ref{lemma3}, we get
$$
\Theta<mf\sqrt{\alpha_0}+\frac{m-1}{\sqrt{\alpha_0}}\delta f+\delta\sqrt{\frac{m-1}{m}}.
$$

Now, we apply inequality \eqref{ineqTerm1} twice. First, one obtains
$$
\delta\sqrt{\frac{m-1}{m}}<B_m f\sqrt{\alpha_0},
$$
and then, it quickly follows that $f<\sqrt{\alpha_0(m-1)/m}$.

Next, using the later inequality and, again, \eqref{ineqTerm1}, one obtains
$$
\frac{m-1}{\sqrt{\alpha_0}}\delta f<(m-1)B_m f\sqrt{\alpha_0}.
$$
Therefore,
\begin{eqnarray*}
% \nonumber % Remove numbering (before each equation)
\Theta &<& mf\sqrt{\alpha_0}+(m-1)B_m f\sqrt{\alpha_0}+B_m f\sqrt{\alpha_0}=  \\
&=& m\left(1+B_m\right)f\sqrt{\alpha_0}.
\end{eqnarray*}
Moreover, since $\alpha_0\leq m_0$ and $1+B_m\leq 6/5$ (from \eqref{Bm}), we get
\begin{equation}\label{ineqTheta}
\Theta<\frac{6}{5}mf\sqrt{m_0}.
\end{equation}

We continue our argument by finding an appropriate upper bound for the last term of \eqref{eq:deltaLambda1}. Taking into account Lemma \ref{lemma2} and the expression of $\left|\nabla A\right|^2$ with respect to $b_{ijk}$'s, it follows that, at $p$,
\begin{eqnarray}\label{ineqSum}
\nonumber \sum\limits_{\substack{i=1 \\ k\geq 2}}^m \frac{b_{i1k}^2}{\mu_1-\mu_k} &\geq & \sum\limits_{\substack{i=1 \\ k\geq 2}}^m \frac{b_{i1k}^2}{\mu_1-\mu_2}> \\
\nonumber &>& -\frac{3-3^{-10}}{2}\sqrt{\frac{m-1}{m}}\sum\limits_{\substack{i=1 \\ k\geq 2}}^m b_{i1k}^2 =\\
\nonumber &=& -\frac{1-3^{-11}}{2}\sqrt{\frac{m-1}{m}}\left(3\sum\limits_{\substack{i=1 \\ k\geq 2}}^m b_{i1k}^2\right) \geq\\
\nonumber &\geq & -\frac{1-3^{-11}}{2}\sqrt{\frac{m-1}{m}}\left(3\sum\limits_{\substack{i=1 \\ k\geq 2}}^m b_{i1k}^2 +b_{111}^2+\sum_{i,j,k\geq 2}^{m}b_{ijk}^2\right)=\\
& = & -\frac{1-3^{-11}}{2}\sqrt{\frac{m-1}{m}}\left|\nabla A\right|^2.
\end{eqnarray}
Now, from \eqref{eq:deltaLambda1}, \eqref{ineqTheta} and \eqref{ineqSum}, we achieve, at $p$,
\begin{equation}\label{ineqDeltaLambda1}
\Delta \lambda_1<\frac{6}{5}m m_0\eta f-3^{-11}\sqrt{\frac{m-1}{m m_0}}\left|\nabla A\right|^2-m_0\phi f\sqrt{m(m-1)}.
\end{equation}
We are not yet able to conclude that $\Delta \lambda_1<0$ on $M$, so we need a better estimation. We can continue our process, with the following algebraic remarks: $\sqrt{m(m-1)}<m+1$, for any positive integer $m$, $3\eta+2\sigma<3/4$, from Lemma \ref{lemma2}, and $9m/4-3/4>2m$, for any $m>3$. Then, using these inequalities and \eqref{eq2Lemma1}, we obtain
\begin{eqnarray*}
% \nonumber % Remove numbering (before each equation)
\phi\sqrt{m(m-1)} &>& \eta\left[3m-3(m+1)\eta-2(m+1)\sigma\right]=\\
&=& \eta\left[3m-(m+1)(3\eta+2\sigma)\right]> \\
&>& \eta\left[3m-\frac{3}{4}(m+1)\right]=\\
&=& \eta\left(\frac{9}{4}m-\frac{3}{4}\right)>\\
&>& 2m\eta,\\
\end{eqnarray*}
and, therefore
$$
\frac{6}{5}m m_0\eta f<\frac{3}{5}m_0\phi f\sqrt{m(m-1)}.
$$
From \eqref{ineqDeltaLambda1} and the above estimation, we conclude with
\begin{equation}\label{ineqDeltaLambda1-2}
\Delta \lambda_1<-3^{-11}\sqrt{\frac{m-1}{m m_0}}\left|\nabla A\right|^2-\frac{2}{5}m_0\phi f\sqrt{m(m-1)}.
\end{equation}
Therefore, at $p$, we have $\Delta\lambda_1<0$. As the point $p$ was arbitrarily fixed, we conclude that $\Delta \lambda_1<0$ on $M$.

Further, since $|A|^2=m$, we obtain that the Ricci curvature of $M$ is bounded below,
$$
\ricci(X,X)\geq-2m(m-1), \qquad X\in C(TM).
$$
Knowing also that $M$ is complete, we can apply the Omori-Yau maximum principle (see, for example \cite{CY75, O67, Y75}) and obtain that there exists a sequence of points $\left\{p_k\right\}_{k\in\mathbb{N}}\subset M$ such that
$$
\left(\Delta\lambda_1\right)\left(p_k\right)>-\frac{1}{k}.
$$
But we have seen that $\Delta \lambda_1<0$ at any point of $M$, so, in particular, $\left(\Delta\lambda_1\right)\left(p_k\right)<0$. Therefore
$$
\lim_{k\to\infty}\left(\Delta\lambda_1\right)\left(p_k\right)=0,
$$
and, moreover, from \eqref{ineqDeltaLambda1-2}, we deduce
$$
\lim_{k\to\infty}|\nabla A|^2\left(p_k\right)=0 \quad \text{ and } \quad \lim_{k\to\infty}\phi\left(p_k\right)=0.
$$
Now, using the fact that the quantity $\left|\nabla A\right|^2+mf\phi\sqrt{m_0^3}$ is constant (see \eqref{eq1}), it follows that
$$
|\nabla A|\equiv 0 \quad \text{ and } \quad \phi\equiv 0.
$$
Further, on the one hand, from \eqref{eq1} we get
$$
m_0-m\left(f^2+1\right)+(m-2)f\sqrt{\frac{m}{m-1}m_0}=0,
$$
and on the other hand, from \eqref{eqLemma3}, we have
$$
\alpha_0-m\left(f^2+1\right)+(m-2)f\sqrt{\frac{m}{m-1}\alpha_0}=0.
$$
Therefore, $\alpha_0=m_0$ and then
$$
f=\frac{m-2}{m}.
$$

Finally, from $\phi\equiv 0$, we see that we have equality in the Okumura Lemma (in the left-hand side), see \cite{O74}, and so
$$
\lambda_1=-1 \qquad \text{and} \qquad \lambda_2=\lambda_3=\cdots=\lambda_m=1.
$$
Thus, $M$ is the Clifford torus $\mathbb{S}^{1}\left(1/\sqrt{2}\right)\times \mathbb{S}^{m-1}\left(1/\sqrt{2}\right)$.

\end{document}